\documentclass[12pt]{article}

\usepackage{amssymb,amsmath,amsfonts,amssymb}
\newcommand{\be}{\begin{equation}}
\newcommand{\ee}{\end{equation}}
\newcommand{\la}{\label}
\newcommand{\ba}{\begin{array}{c}}
\newcommand{\ea}{\end{array}}
\newcommand{\Rm}{{\mathbb R}}
\renewcommand{\d}{{\mbox{div}}_g}
\newcommand{\D}{ {\Delta_g}}
\newcommand{\nag}{\nabla_g}
\newcommand{\nax}{\nabla_x}
\newcommand{\Dx}{\Delta_x}
\newcommand{\dx}{{\mbox{div}}_x}

\newtheorem{thm}{Theorem}
\newtheorem{prop}{Proposition}
\newtheorem{lemma}{Lemma}
\title{Regularity of coupled two-dimensional\\Nonlinear Fokker-Planck and Navier-Stokes Systems}
\author{P. Constantin\\Department of Mathematics, The University of Chicago\\ Chicago IL 60637,\\C. Fefferman \\Department of Mathematics Princeton University\\ Princeton NJ 08544,\\ E.S. Titi\\Department of Mathematics\\Department of Mechanics and Aerospace Engineering \\Univeristy of California\\Irvine CA 92697  and \\Department of Computer Science and \\ Applied Mathematics\\Weizmann Institute of Science \\ Rehovot 76100 Israel,\\ A. Zarnescu
\\Department of Mathematics\\The University of Chicago\\Chicago, IL 60637}

\begin{document}

%version{04/16/06}

\maketitle
\newpage
\noindent{\bf Abstract.} We consider systems of particles coupled with fluids. The  particles are described by the evolution of their density, and the fluid is described by the Navier-Stokes equations. The  particles add stress to the fluid and the fluid carries and deforms the particles. Because the particles perform rapid random motion, we assume that the density of particles is carried by a time average of the fluid velocity. The resulting coupled system is shown to have smooth solutions at all values of parameters, in two spatial dimensions.



\section{Introduction}
We discuss global regularity of solutions of systems of equations describing
fluids with particle suspensions. The particles are parameterized by
independent microscopic variables $m$  that belong to a compact,
connected, smooth Riemannian manifold $M$ of dimension $d$. Derivatives with
respect to the microscopic variables are designated by the subscript
$g$. The particles are included in a fluid in ${\mathbb{R}}^n,\,\, n=2$
obeying the forced Navier-Stokes equations.

The forces exerted by the particles on the fluid are expressed
through the divergence of an added stress tensor. The added stress
tensor $\tau_p(x,t)$ is obtained after averaging out the microscopic variable and  the Navier-Stokes equation is macroscopic. The
microscopic inclusions at time $t$ and macroscopic physical location
 $x\in {\mathbb R}^2$ are described by the  density  $f(x,m,t)dm$
where $dm$ is the Riemannian volume element in $M$. The density is nonnegative, $f\ge 0$ and 
\be \rho(x,t) = \int_M f(x,m,t)dm \le 1 \la{intone} \ee 
holds for every $x$, $t\ge 0$.

The added stress tensor is given by an expansion 
\be \tau_p(x,t) =
\sum_{k=1}^{\infty}\tau^{(k)}(x,t) \la{expa} \ee 
where 
\be
\tau^{(1)}_{ij}(x,t) = \int_M \gamma_{ij}^{(1)}(m) f(x,m,t)dm,
\la{sigmaone} 
\ee 
\be \tau^{(2)}_{ij}(x,t) =
\int_M\int_M\gamma_{ij}^{(2)}(m_1,m_2)f(x,m_1,t)f(x,m_2,t)dm_1dm_2,
\la{sigmatwo} \ee and, 
in general, 
\be \ba \tau_{ij}^{(k)}(x,t) =
\\  \int\limits_{M\times\cdots\times M}\gamma_{ij}^{(k)}(m_1,\dots,
m_k)f(x, m_1, t)f(x,m_2,t)\cdots f(x, m_k, t)dm \ea
\la{sigmak} 
\ee 
Expansions of this kind for the added stress tensor $\tau_p$ are encountered
in the polymer literature (\cite{doied}). In (\cite{c}) it was
proved that only two  structure coefficients in the expansion,
$\gamma^{(1)}_{ij}, \gamma^{(2)}_{ij}$ are needed in order to have
energetically balanced equations, provided certain constitutive
relations are imposed. The energy balance confers stability to
certain time-independent solutions of the equations. In this work we
are interested only in general existence results, and do not need to
use special constitutive relations. We will only use the fact that
the coefficients $ \gamma^{(k)}_{ij}$ are
 smooth, time independent, $x$ independent, $f$ independent. When
infinitely many coefficients are present, we will use a finiteness
condition assuming that the series 
\be
\sum_{k=1}^{\infty}k^3\|\gamma^{(k)}_{ij}\|_{H^{\rho_k}(M\times\cdots\times
M)} \la{sigmasum} \ee 
converges for a sequence  $\rho_k>  \frac{k+
4d + 6}{2}$.

From (\ref{intone}, \ref{expa}, \ref{sigmak}, \ref{sigmasum} ) 
it follows that 
\be |\tau_p(x,t)|\le c\rho(x,t) \la{sigmalinfty} \ee 
holds with a constant that depends only
on the coefficients $\gamma^{(k)}_{ij}$. The spatial gradients of $\tau_p$
are of particular importance for regularity. The fact that the constitutive
coefficients $\gamma^{(k)}_{ij}$ are smooth functions of the microscopic variables allows us to relate the size of the spatial gradients of $\tau_p$ to a rather coarse average on $M$: differentiating
(\ref{sigmak}) with respect to $x$ it follows from (\ref{sigmasum})
that 
\be |\nax\tau_p (x,t)| \le cN(x,t)\la{naxsigma} \ee 
holds with
a constant $c$ that depends only on the smooth coefficients 
$\gamma^{(k)}_{ij}$. Here
$$
N(x,t) = \|R\nax f\|_{L^2(M)}
$$
and
$$
R = \left (-\Delta_g + {\mathbb I}\right)^{-\frac{s}{2}}
$$
with $s>\frac{d}{2} +1$. 
The inequalities (\ref{sigmalinfty}) and (\ref{naxsigma}) are the
only information concerning the relationship between $\tau_p$ and
$f$ that we need for regularity results in this
work. We can use the detailed form (\ref{expa}) and the finiteness
condition (\ref{sigmasum}) to deduce them, but we could just as well
require them instead of (\ref{expa}).

The particles are carried by the fluid, agitated by thermal noise and interact among themselves in a mean-field fashion, through potentials
that depend linearly and nonlocally on the particle density
distribution $f$ (\cite{ons}). We assume that the fluid
does not vary much in time during a characteristic relaxation time of  the particles. Mathematically, this means that the particles are carried by a 
short time average of the fluid velocity. This assumption alllows us to
prove global existence of smooth solutions and to bound
a priori the size of the physical space gradients of the stresses.

The mathematical study of complex fluids is in a developing stage.
Most results are for models that are macroscopic closures, that is, 
in which $\tau_p$ has its own macroscopic evolution, coupled with 
the fluid: the microscopic variables do not appear at all.
Existence theory for viscoelastic Oldroyd models is presented in (\cite{linliuzang}); see also (\cite{sideris}) for related issues.  
There are few other regularity results concerning complex fluids, including some that retain microscopic variables. Among them are (\cite{ez}, \cite{jourdain}, \cite{jourdainl}, \cite{lizz}, \cite{renardy}). For Smoluchowski equations coupled with fluids,
the case in which $u$ is given by a time independent  linear Stokes equation in $n=3$, $M={\mathbb S}^2$ with $\tau_p$ given by a relation (\ref{sigmaone}) was studied in (\cite{otto-tzavaras}) for the case of a linear Fokker-Planck equation and in (\cite{c}) for general nonlinear Fokker-Planck equations.

The proofs in the present work are based on a few key facts. The first one is that gradients of $\tau_p$ are bounded by $N$, and $N$ is controlled linearly by the advecting velocity, taking advantage of the a priori boundedness of $\rho$ in $L^1\cap L^{\infty}$. This idea was used in (\cite{c}) to prove regularity for the system of particles coupled with the Stokes system for the fluid. The second significant fact concerns the Navier-Stokes system driven by the divergence of  bounded stresses. We are interested in the size of the time integral of the supremum of the norm of the gradient of velocity. This is an important nondimensional magnitude that controls the amplification of gradients of passively advected scalars. We prove a logarithmic bound for this amplification factor. The strategy of proof uses a natural idea 
introduced in (\cite{CHMA01}): time integration is performed first in  each wave-number shell, to take advantege of the rapid smoothing of small scales due to viscosity.

\section{Estimates for  2-D Navier-Stokes Equations}

\par Consider a Navier-Stokes system in ${\mathbb R}^2$:

\begin{eqnarray}\setlength\arraycolsep{2pt}
\frac{\partial u}{\partial t}+u\cdot\nabla u- \Delta u + \nabla p&=&\nabla\cdot \sigma,\nonumber\\
\nabla\cdot u&=&0, \nonumber\\
u|_{t=0}&=&u_0, \label{eq:ns}
\end{eqnarray} 
where $\sigma (x,t)$ is a symmetric two-by-two matrix that, in this section,  will be considered to be a given function. We will be interested in estimates when $\sigma$ is integrable and bounded
by constants that are known a priori and are of order one. The physical space gradients of $\sigma$ are possibly large. The aim of the bounds is to find the effect that these gradients have on the stretching amplification term
$$
\int_0^t\|\nabla u (t)\|_{L^{\infty}}dt.
$$
We take smooth, divergence-free initial velocities
$$
u(0)\in L^{2}({\mathbb R}^2)\cap W^{2,r}({\mathbb R}^2) 
$$
with localized vorticity $\omega= \nabla^{\perp}\cdot u$,
$$
\omega(0)\in L^2({\mathbb R}^2)\cap W^{1,r}({\mathbb R}^2)
$$
with $r>2$.
We recall the energy estimate
\be
\sup_{t\le T}\|u(t)\|_{L^2(dx)}^2 + \int_0^T\|\nabla u(t)\|_{L^2(dx)}^2dt \le
\int_0^T\|\sigma(t)\|_{L^2(dx)}^2dt + \|u(0)\|_{L^2}^2\la{kine}
\ee
and the fact that in two dimensions the vorticity obeys
\be
\partial_t \omega + u\cdot\nax \omega - \Delta\omega =
\nax^{\perp}\cdot\dx\sigma.\la{vorteq} 
\ee 
\begin{lemma} Let $r\ge 2$. There exists a constant $c_r$  such that
\be
\sup_{t\le T}\|\nabla u(t)\|_{L^r(dx)}^2 \le c_r\int_0^T\|\nabla \sigma(t)\|_{L^r}^2dt + c_r\|\omega(0)\|_{L^r}^2
\la{grr}
\ee
holds.
For $r=2$ we have additionally,
\be
\int_0^T\|\Delta u(t)\|^2dt \le c_2\int_0^T\|\nabla\sigma\|_{L^2(dx)}^2 + \|\omega(0)\|_{L^2}^2.
\la{lap}
\ee

\end{lemma}

\par\noindent{\bf Proof.}
We multiply
(\ref{vorteq}) by $\omega^{r-1}$ , $r\ge 2$ 
and integrate by parts to obtain
$$
\frac{d}{rdt}\|\omega \|_{L^r(dx)}^r + \int\omega^{r-2}|\nax
\omega|^2dx \le (r-1)\int|\nax \sigma||\nax\omega|\omega^{r-2}
$$
We use a H\"{o}lder inequality with exponents $r, 2, \frac{2r}{r-2}$,
and then with exponents $2,2$:
\begin{eqnarray}
\frac{d}{r dt}\|\omega \|_{L^r(dx)}^r + \int\omega^{r-2}|\nax
\omega|^2dx \nonumber\\
\le (r-1)\|\nax\cdot\sigma\|_{L^r}\left(\int\omega^{r-2}|\nax
\omega|^2dx\right)^{1/2} \|\omega\|_{L^r}^{\frac{r-2}{2}}dx
\la{zh}\nonumber\\
\le
\frac{(r-1)^2}{2}\|\nax\cdot\sigma\|_{L^r}^2\|\omega\|_{L^r}^{r-2}+
\frac{1}{2}(\int\omega^{r-2}|\nax\omega|^2dx)
\end{eqnarray} which implies that
\be \sup_{t\in [0,T]} \|\omega(t)\|_{L^r}^2\le
{C(r-1)^2}\int_0^T
\|\nax\cdot\sigma(t)\|_{L^r}^2dt + \|\omega(0)\|_{L^r}^2\la{omegalr}\ee
and thus (\ref{grr}) follows. In order to obtain (\ref{lap}) we integrate (\ref{zh}) in time at  $r=2$.

\par We need  a logarithmic inequality  for $\|u\|_{L^{\infty}}$. Such inequalities were first introduced in
(\cite{brega}). We will write $\log_{*}(\lambda) = \log(2+\lambda)$ for $\lambda>0$. Note that $\log_{*}(0)>0$ and $\log_{*}(\lambda\mu)\le
\log_{*}(\lambda) + \log_{*}(\mu)$ holds for $\lambda\ge 0, \,\mu\ge 0$. We check our inequality   
\be
\|u(t)\|_{L^{\infty}}\le
C_r\|\omega(t)\|_{L^2}\left \{1+\sqrt{\log_{*}\left\{ \left(\frac{\|\omega(t)\|_{L^r}}{\|\omega(t)\|_{L^2}} \right )^{\frac{r}{r-2}}\frac{\|u(t)\|_{L^2}}{\|\omega(t)\|_{L^2}}\right\}}\right \}
\la{bg}
\ee
directly from the Biot-Savart law:
$$
u(x,t) = \frac{1}{2\pi}\int_{{\mathbb R}^2} \left ({\widehat {z}}^{\perp}\right )\omega(x-z,t) \frac{dz}{|z|}
$$
where ${\widehat{z}} = \frac{z}{|z|}$.
We pick two numbers $0<l\le L$, take a smooth radial function $a(s)$ , $0\le a(s)\le 1$,  that equals identically one for $s\le 1$ and identically zero for $s\ge 2$, and write 
$$
u^{(l)}(x,t) =\frac{1}{2\pi} \int_{|z|\le l}\left ({\widehat {z}}^{\perp}\right )\omega(x-z, t) \frac{dz}{|z|},
$$
$$
u_{(l)}^{(L)}(x,t) =\frac{1}{2\pi} \int_{|z|\ge l}\left ({\widehat {z}}^{\perp}\right )\omega(x-z, t) a\left(\frac{|z|}{L}\right)\frac{dz}{|z|}
$$
and
$$
u_{(L)}(x,t) =\frac{1}{2\pi} \int_{|z|\ge l}\left ({\widehat {z}}^{\perp}\right )\omega(x-z,t) \left (1-a\left (\frac{|z|}{L}\right)\right)\frac{dz}{|z|}
$$
Clearly
$$
u = u^{(l)} + u_{(l)}^{(L)} + u_{(L)}
$$
holds pointwise. It is also clear that
$$
|u^{(l)}(x,t)| \le \|\omega(t)\|_{L^r}l^{\frac{r-2}{r}}
$$
and that
$$
|u_{(l)}^{(L)}(x,t)| \le \|\omega(t)\|_{L^2} \log_{*}\left (\frac{2L}{l}\right).
$$
We integrate by parts in the term $u_{(L)}$, using $\omega = \nabla^{\perp}\cdot u$ and deduce
$$
|u_{(L)}(x,t)| \le C\frac{1}{L}\|u(t)\|_{L^2}
$$
We choose 
$$
l= \left(\frac{\|\omega(t)\|_{L^2}}{\|\omega(t)\|_{L^r}}\right)^{\frac{r}{r-2}}
$$
and
$$
L = \frac{\|u(t)\|_{L^2}}{\|\omega(t)\|_{L^2}}
$$
if, with this choice, it turns out that $l<L$. If not, then we still take $L$ as above, but we take $l=L$. The inequality (\ref{bg}) follows. 
\par Let us consider
\be 
K_r^{(p)} = \|\sigma\|_{L^p(0,T;L^r)} = \left\{\int_0^T\|\sigma(t)\|_{L^r(dx)}^pdt\right\}^{\frac{1}{p}} 
\la{krp}
\ee
and 
\be
B_{r}^{(p)} = \|\nabla\cdot\sigma\|_{L^p(0,T; L^r(dx))} = \left\{\int_0^T\|\nabla\cdot \sigma (t)\|_{L^r}^pdt\right\}^{\frac{1}{p}}.
\la{br}
\ee
Note that (\ref{kine}) implies that
\be
\sup_{t\le T}\|u(t)\|_{L^2}^2 + \int_0^T\left \|\nabla u(t)\right\|_{L^2}^2dt \le K_0.
\la{intgru}
\ee
with 
\be
K_0 = (K_2^{(2)})^2 + \|u_0\|_{L^2}^2\la{ebound}
\ee
and note that (\ref{grr}) implies that
\be
\sup_{t\le T}\|\omega(t)\|_{L^r} \le c_r(B_r^{(2)} + \|\omega(0)\|_{L^r}) = \Omega_r.\la{bgrr}
\ee
Finally, using (\ref{lap}) and (\ref{intgru}) we can write
\be
\int_0^T\|u(t)\|^2_{H^2}dt\le c(1+T)\left\{ \|\sigma\|^2_{L^2(0,T; H^1)} + \|u(0)\|^2_{H^1}\right\}.
\la{htwo}
\ee

Now we integrate the square of (\ref{bg}) in time, taking the supremum in time of the logarithmic part using (\ref{bgrr}), and bounding the time integral of the sqaure of the gradients using (\ref{intgru}).
\begin{lemma} For $r>2$ there exists a constant $c_r$ such that
\be
\int_0^T\|u(t)\|^2_{L^{\infty}}dt \le c_rK_0\left\{\log_{*}(\Omega_r) + \log_{*}(\Omega_2) +\log_{*} K_0\right\}\la{unftyb}
\ee
\end{lemma}

\begin{thm} For $r>2$ there exists a constant $c_r$ such that
\be
\ba
\|\nax\nax u\|_{L^2(0,T; (L^r))}\le \\
c_rT\|\nabla\omega_0\|_{L^r} + c_r{\sqrt{K_0}}\Omega_r\sqrt{\left\{\log_{*}(\Omega_r) + \log_{*}(\Omega_2) +
\log_{*}(K_0)\right\}}
\ea
\la{grgrr}
\ee
holds.
\end{thm}
\noindent{\bf Proof.}
We represent
\be
\nax\omega (t) = {\mathcal T_1}(\nax^{\perp}\cdot\sigma) - 
{\mathcal T_2}(u \omega ) + e^{t\Delta}(\nax\omega(0))
\la{naxomega}
\ee
where ${\mathcal T_1}$ and ${\mathcal T_2}$ are operators of the form
$$
h(t) \mapsto {\mathcal T}h = \int_0^t e^{(t-s)\Delta}\Delta
{\mathbb H}h(s)ds
$$
with ${\mathbb H} = {\mathbb H}(D)$ homogeneous of degree zero. Such
operators are bounded in $L^p(dt; L^q(dx))$ for $1<p, q<\infty$ by the maximal regularity of the heat equation combined with the boundedness of the operators
${\mathbb H}$ in $L^q$ spaces (see, for example (\cite{lemarie})).
The inequality follows then from the bound
\be
\int_0^T \|u(t)\omega(t)\|^2_{L^r}dt \le c_r K_0\Omega_r^2\left\{\log_{*}(\Omega_r) + \log_{*}(\Omega_2) + \log_{*}(K_0)\right\}
\la{uomegab}
\ee
which, in turn, follows from (\ref{bgrr}) and (\ref{unftyb}).

\smallskip
We mention also an a priori bound for $\int_0^T\|u\|^p_{L^{\infty}}dt$
for $p<2$.
This is obtained as follows: we write
\be
\ba
\omega(t) = e^{t\Delta}\omega_0 + \\
+ \int\limits_0^t e^{(t-s)\Delta}\left\{\nabla^{\perp}\cdot\dx\sigma(s) +\partial_2\partial_1( 
u_2^2 -u_1^2)(s) + (\partial_2^2-\partial_1^2)u_1u_2(s)\right\}ds\la{omegu}
\ea 
\ee
This follows from a well-known identity 
$$
u\cdot\nabla\omega =  -\partial_2\partial_1( 
u_2^2 -u_1^2) - (\partial_2^2-\partial_1^2)u_1u_2
$$
An easy calculation verifies this after one writes $u_1 = -\partial_2\psi$, $u_2
= \partial_1\psi$, $\omega = \Delta \psi$.
From (\ref{omegu}) we get, for any $1< p,\,q <\infty$
\be
\|\omega\|_{L^p[(0,t); L^q(dx)]}\le C_{pq}(t) + C_{pq}\|\sigma\|_{L^{p}[(0,t); L^q(dx)]} +
C_{pq}\|u\|^2_{L^{2p}[(0,t), L^{2q}(dx)]}
\la{oepq}
\ee
We know however from (\ref{intgru}) that
$$
\|u\|_{L^{\infty}[(0,t);L^2(dx)]} \le K_0
$$
and
$$
\|u\|_{L^2[(0,t); L^r(dx)]} \le C_r
$$
a priori, for any $r<\infty$, with $K_0, C_r$ independent of $t\in [0,T]$. (These constants may depend on $T$ because the norm of $\sigma$ in $L^2$ may depend on $T$).
Then, by interpolation
\be
\|u\|_{L^{p}[(0,t); L^q(dx)]} \le C_{pq}
\la{quineq}
\ee
holds for $q\ge 2$ and $p<\frac{2q}{q-2}$. In view of \ref{oepq}
we get that
\be
\|\omega\|_{L^p[(0,t); L^q(dx)]}\le C_{pq}
\la{ompq}
\ee
holds for $q\ge 2$ and $p< \frac{q}{q-1}$. Then, taking $q>2$ and using a Sobolev embedding theorem, we obtain that
\be
\int_0^T \|u\|_{L^{\infty}(dx)}^pdt \le C_p
\la{almost}
\ee
holds a priori, for any $p<2$.

\smallskip

For the bound for $\int_0^T\|\nabla u(t)\|_{L^{\infty}}dt$ 
we will use the Littlewood-Paley decomposition. Let $\mathcal{D}(\Omega)$ denote the set of
$C^\infty$ functions compactly supported in
$\Omega$. Let  $\mathcal{C}$ be the annulus centered at
$0$, and with radii $1/2$ and $2$. There exist two
nonnegative, radial functions $\chi$ and $\varphi$, belonging
respectively to $\mathcal{D}(B(0,1))$ and to
$\mathcal{D}(\mathcal{C})$ so that

\begin{displaymath}
\chi(\xi)+\Sigma_{j\ge 0}\varphi(2^{-j}\xi)=1, 
\end{displaymath}
and
\begin{displaymath}
|j-k|\ge 2\Rightarrow \sup (\varphi(2^{-j}))\cap\sup
(\varphi(2^{-k}))=\varnothing.
\end{displaymath}

\par We denote by $\mathcal{F}$ the Fourier transform on
$\mathbb{R}^2$ and  let $h,\tilde h,\Delta_j,S_j(j\in {\mathbb N})$ be
defined by

\begin{displaymath}
h=\mathcal{F}^{-1}\varphi\,\,\,\textrm{and}\,\,\,\tilde
h=\mathcal{F}^{-1}\chi,
\end{displaymath}

\begin{displaymath}
\Delta_j u=\mathcal{F}^{-1}(\varphi(2^{-j}\xi)\mathcal{F}u)=2^{2j}
\int h(2^jy)u(x-y)dy,
\end{displaymath}

\begin{displaymath}
S_j u=\mathcal{F}^{-1}(\chi(2^{-j}\xi)\mathcal{F}u)=2^{2j}
\int\tilde h(2^j y)u(x-y)dy.
\end{displaymath}

\par Then

\begin{displaymath}
u=S_0 u + \sum_{j\ge 0}\Delta_j(u)
\end{displaymath} where $u\in\mathcal{S}'$, the space of tempered
distributions, and the equality holds in the sense of distributions.

\smallskip\par The well-known Bernstein inequalities
(see, for instance \cite{CH98}) express the fact that $\Delta_j$ is localized
around the frequency $2^j$.

\begin{prop}(Bernstein Inequalities) Let $a,b\in [1,\infty]$, $j\ge 0$.
There exists constants $c$, independent of $a,b, j$ such that the following hold:

\begin{equation}
\|\Delta_j e^{t\Delta}u\|_{L^a}\le ce^{-t2^{2{(j-1)}}}\|\Delta_j
u\|_{L^a}, \label{eq:localizedheat}
\end{equation}

\begin{equation}
\|\Delta_j f\|_{L^\infty}\le c\|f\|_{L^\infty},
\label{eq:localizedlinfty}
\end{equation}
and, in two space dimensions:
\be
\|S_j f\|_{L^{\infty}} \le c \left (\|f\|_{L^2} + \sqrt{j}\|\nabla f\|_{L^2}\right)\la{sqb}
\ee
and

\begin{equation}
\|\Delta_j \partial^\alpha u\|_{L^a}\le
c2^{j|\alpha|+2j(1/b-1/a)}\|\Delta_j u\|_{L^b}
\label{eq:localizedderivative}
\end{equation} where $|\alpha|$ is the length of the
multiindex $\alpha$. \label{lemma:bernstein}
\end{prop}
The Littlewood-Paley decomposition is best suited for Besov spaces $B^{s}_{p,q}$ defined by requiring the sequence $2^{sj}\|\Delta_j(u)\|_{L^p}$ to belong to $\ell_q$ and by requiring $S_0(u)$ to be in $L^p$.  
$L^2$ based Sobolev space norms can be  computed in terms of the Littlewood-Paley
decomposition:
$$
\|u\|_{H^s}^2 \sim \|S_0(u)\|_{L^2}^2 + \sum_{j\ge 0}2^{2sj}\|\Delta_j u\|^2_{L^2}
$$
where $\sim$ means equivalence of norms. 
However, the norm we are interested in is the $L^1(0,T; W^{1,\infty})$ norm.
The $C^s$ norms can be computed as
$$
\|u\|_{C^s}\sim \|S_0(u)\|_{C^s} + \sup_{j\ge 0}2^{js}\|\Delta_j (u)\|_{L^{\infty}}
$$ 
but only if $s$ is not an integer. In order to obtain $L^{\infty}$ bounds
for the gradient we will have to resort to the inequality:
$$
\|\nabla u\|_{L^{\infty}} \le \|S_0(\nabla u)\|_{L^{\infty}} + \sum_{j\ge 0}\|\Delta_j(\nabla u)\|_{L^{\infty}}.
$$
This inequality reflects the embedding $B^{0}_{\infty, 1}\subset L^{\infty}$, which is a strict inclusion. The advantage of using this sum (the norm in $B^{0}_{\infty,1}$) is that we can commute time integration and summation, while time integration and supremum do not commute in general.

\medskip
\par\noindent 

\begin{thm}
Let $u$ be  a solution of the $2D$ Navier-Stokes system
(\ref{eq:ns}),
 with divergence-free initial data $u_0\in W^{1,2}({\mathbb R}^2)\cap W^{1,r}({\mathbb R}^2)$. Let $T>0$ and let the forces $\nabla\cdot\sigma$ obey
$$
\sigma\in L^1(0,T;L^{\infty}({\mathbb R}^2))\cap L^{2}(0,T; L^2({\mathbb R}^2))
$$
and
$$
\nabla\cdot\sigma \in L^1(0,T; L^r({\mathbb R}^2))\cap L^2(0,T; L^2({\mathbb R}^2))
$$
with $r>2$. There exists a constant $c$ depending on $r$ such that, for every
$\epsilon>0$
\be
\ba
\int_0^T \|\nabla u\|_{L^\infty}dt\le c\sqrt{T}\|u(0)\|_{H^1} + \\
 c\left\{K_2^{(1)} T +
K_{\infty}^{(1)}\log_{*}\left (\frac{B_r^{(1)}}{\epsilon}\right) +(1+T)E_1\log_{*}\left (\frac{(1+T)R_1}{\epsilon }\right) + \epsilon\right\}
\ea
\la{gineq}
\ee
holds, where 
\be
E_1 = \int_0^T \|u(t)\|_{H^1}^2dt.
\la{eone}
\ee
and
\be
R_1 = \left\{\|\sigma\|_{L^2(0,T; H^1)}^2 + \|u(0)\|^2_{H^1}\right\}.
\la{Rone}
\ee 
Consequently, in view of (\ref{intgru}) and (\ref{gineq}) above
\be
\int_0^T\|\nabla u(t)\|_{L^{\infty}}dt \le c(1+T)^2K\log_{*}(B)
\la{vague}
\ee
with $K= K_0+K_{\infty}^{(1)} + K_{2}^{(1)} + \|u(0)\|_{H^1}$ depending on norms of $\sigma$ and the initial velocity, but not on gradients of $\sigma$, and $B= B_r^{(1)} + 
R_1$ depending on norms of the gradients of $\sigma$.
\end{thm}

\smallskip\noindent{\bf Remark.}  
The bound is in fact for the stronger norm of $u$ in the inhomogeneous space $L^1(0,T; B^1_{\infty, 1})$.

\smallskip\noindent{\bf Proof.} We start with the Duhamel formula
 for the gradient of solutions of (\ref{eq:ns})

\be
 \ba
  \nabla u =e^{t\Delta }\nabla u_0 + \int_0^t e^{(t-s)\Delta}\Delta{\mathbb H}(D)({u(s)\otimes u(s))}ds\\
  -\int_0^t e^{(t-s)\Delta}\Delta {\mathbb{H}(D)}(\sigma(s))ds
\ea
\label{nablau}
\ee
The homogeneous operator ${\mathbb H}(D)$ 
is given by
\be
{\mathbb H}(D)(u\otimes u)_{ij} = R_j(\delta_{il} + R_iR_l)R_k(u_lu_k)
\la{hs}
\ee
with $R_j = \partial_j(-\Delta)^{-\frac{1}{2}}$ Riesz transforms.
The strategy is based on an idea of Chemin and Masmoudi (\cite{CHMA01}) to take the time integral first, for each wave number shell. They did not use information about derivatives of $\sigma$, and therefore obtained only bounds for $\sup_{q\ge 1}\int_0^T\|\Delta_q \nabla u\|_{L^{\infty}}ds$. We will use the gradients of $\sigma$ to bound the high frequencies and will sum in $q$ in order to
estimate $\|\nabla u\|_{L^1((0,t); L^{\infty})}$. Also, we use somewhat 
different estimates than them for the individual shell contributions, but 
like them, we take advantage of a time integration at each shell.
We treat separately the contributions coming from $\sigma$ and those coming from $u\otimes u$:
$$
\nabla u = F + U + e^{t\Delta}\nabla u_0
$$
where 
\begin{equation}
 F(t) = -\int_0^t e^{(t-s)\Delta}\Delta {\mathbb H}(D)(\sigma(s))ds.
 \label{eq:F}
 \end{equation}
and
\be
U(t) = \int_0^t e^{(t-s)\Delta}\Delta {\mathbb H}(D)((u\otimes u)(s))ds.
 \la{Ut}
 \end{equation}
Clearly
\be
\int_0^T\|S_0F(t)\|_{L^{\infty}}dt \le cK_2^{(1)} T
\la{szf}
\ee
is true using for instance $\|S_0(F(t))\|_{L^{\infty}(dx)}\le \|(I-\Delta)S_0(F(t))\|_{L^2(dx)}$.
We take $q\ge 0$ and apply $\Delta_q$:

\begin{equation}
\|\Delta_q F(t)\|_{L^\infty}\le \int_0^t
e^{-(t-s)2^{2(q-1)}}2^{2q}\|\sigma(s)\|_{L^\infty}ds\la{dqf}
\end{equation}

\par Integrating on $[0,T]$ and changing order of integration
we obtain
\begin{eqnarray}
\int_0^T \|\Delta_q F(t)\|_{L^\infty}dt\le \int_0^T
\|\sigma(s)\|_{L^\infty}(\int_s^T
e^{-(t-s)2^{2(q-1)}}2^{2q}dt)ds\nonumber\\
\le {c}\int_0^T
\|\sigma(s)\|_{L^\infty}ds. \label{lowF}
\end{eqnarray}

\par We bound the same quantity differently, with large $q$ in mind, and
use $\nabla\sigma\in L^r$ with $r>2$:

\be
\|\Delta_q F(t)\|_{L^\infty} \le c\int_0^t e^{-2^{2(q-1)}(t-s)}
2^{q(1+2/r)}\|\nabla\cdot\sigma\|_{L^r}ds\la{hdqf}
\ee 
We integrate on $[0,T]$ and change the  order of integration as
above, to obtain

\begin{eqnarray}
\int_0^T \|\Delta_q F(t)\|_{L^\infty}dt\le \int_0^T
\|\nabla\cdot\sigma(s)\|_{L^r}(\int_s^T
e^{-2^{2(q-1)}(t-s)}2^{q(1+2/r)}dt)ds\nonumber\\
\le \frac{c}{2^{q(1-2/r)}}\int_0^T \|\nabla\cdot
\sigma(s)\|_{L^r}ds \label{highF}
\end{eqnarray}

Using (\ref{szf}),  (\ref{lowF}) to estimate the small wave numbers in $F$
and (\ref{highF}) to estimate the high ones, we obtain

\be
\ba
\int_0^T\|F(t)\|_{L^{\infty}}dt\le \\
\int_0^T\left [\|S_0 F(t)\|_{L^\infty} +\sum_{0\le q\le M} \|\Delta_q
F(t)\|_{L^\infty}+\sum_{q> M}
\|\Delta_q F(t)\|_{L^\infty}\right ]dt \\
\le cK_2^{(1)}T + cM \int_0^T
\|\sigma(s)\|_{L^\infty}ds + \frac{c} {
2^{M(1/2-1/r)}}\int_0^T \|\nabla\cdot
\sigma(s)\|_{L^r}ds.\ea
\la{fint}
\ee
Then, choosing $M$
\begin{eqnarray}
M = c_r\log_{*} \left(\frac {c\int_0^T \|\nabla\cdot
\sigma(s)\|_{L^r}ds }{\epsilon}\right)
\end{eqnarray}
 we obtain

\be
\int_0^T\|F(t)\|_{L^{\infty}}dt\le  cK_2^{(1)}T + cK_{\infty}^{(1)}\log_{*}
\left(\frac{B_r^{(1)}}{\epsilon}\right) + \epsilon\label{fbound}
\ee
with $K_{\infty}^{(1)}$ defined in (\ref{krp}) and $B_r^{(1)}$ defined in (\ref{br}).

\noindent{\bf Remark.} We do not need to integrate the term $F$ in time, if $\sigma $ is bounded: we can obtain a pointwise logarithmic bound for $F(t)$ in terms of $B_r^{(p)}$ with $p>\frac{2r}{r-2}$.
Indeed, from (\ref{dqf}) we have
\be
\|\Delta_q F(t)\|_{L^{\infty}} \le c \|\sigma\|_{L^{\infty}(dtdx)}
\la{ndqf}
\ee
and from (\ref{hdqf}) we obtain
\be
\|\Delta_q F(t)\|_{L^{\infty}}\le c 2^{-2q(\frac{r-2}{2r}-\frac{1}{p})}\|\nabla\cdot\sigma\|_{L^p(0,T; L^r)}\la{thdqf}
\ee
Summing (\ref{ndqf}) from $q=0$ to $q=M$, summing (\ref{thdqf}) from $q=M$ 
to infinity and choosing $M$ appropriately, we obtain
\be
\sup_{t\le T}\|F(t)\|_{L^{\infty}} \le C\sqrt{K_0T} + cK_{\infty,\infty}\log_{*}{\left (\frac{B_r^{(p)}}{K_{\infty,\infty}} \right )}\la{pointf}
\ee
with
\be
K_{\infty,\infty} = \|\sigma\|_{L^{\infty}(0,T; L^{\infty})}.\la{kii}
\ee
This bound can be used to reprove the global existence of the Stokes system coupled with nonlinear Fokker Planck equations.
 
We split the nonlinear term
\par\noindent  
\be
U(t) = S_2(U)(t) + V(t) 
\la{para}
\ee
with
\be
S_2(U) = S_{0}(U) + \Delta_1(U) + \Delta_2(U),
\la{stwo}
\ee
and 
\be
V(t) = \sum_{q\ge 3}\Delta_q(U)
\la{vt}
\ee
Clearly
\be
\int_0^T\|S_2(U)\|_{L^{\infty}}dt \le c TE_1,\la{s2bound}
\ee
where $E_1$ is given in (\ref{eone}).
For the nonlinear term $V$ we use Bony's decomposition
(see for instance \cite{CH98}) into commensurate and incommensurate
frequencies:
\be
V(t) = C(t) + I(t)\la{vtci}
\ee
with
\be
C(t) = \sum_{q\ge 3}\int_0^te^{(t-s)\Delta}\Delta{\mathbb H}(D)\Delta_q\left(\sum_{|p-p'|\le 2}\Delta_p(u(s))\otimes\Delta_{p'}(u(s))\right)ds
\la{ct}
\ee
and
\be
I(t) = \sum_{q\ge 3}\int_0^te^{(t-s)\Delta}\Delta{\mathbb H}(D)\Delta_q\left (\sum_{|p-p'|\ge 3}\Delta_p(u)\otimes\Delta_{p'}(u)\right)ds
\la{it}
\ee
In the decomposition above we made the convention that the indices $p,p'$ run from $-1$ to infinity and when $p=-1$ instead of $\Delta_p$ we have $S_0$, and of course, the same thing for $p'$.

We start by treating the term $C(t)$. Because the range of $q$ is $q\ge 3$
it follows that the range of $p,p'$ is $p\ge 1$, $p'\ge 1$. Then we
estimate inside the integral 
$$
\ba
\|\Delta_q((\Delta_pu(s))\otimes(\Delta_{p'}(u(s)))\|_{L^{\infty}} \le
c2^{2q}\|\Delta_q((\Delta_pu(s))\otimes(\Delta_{p'}(u(s)))\|_{L^{1}}\\
\le c2^{(2q-2p)}\|\nabla\Delta_pu(s)\|_{L^2}\|\nabla \Delta_{p'}u(s)\|_{L^2}
\ea
$$ 
Clearly at least one of $p,p'$, say $p$, satisfies $p\ge q-2$. 
Now, because $|p-p'|\le 2$, it follows that $p' = p +j$, $j\in [-2, 2]$ and therefore, changing the order of summation in (\ref{ct})
$$
\ba
\|C(t)\|_{L^{\infty}} \le \sum_{j=-2}^2\sum_{p=1}^{\infty}c\int_0^t\|\Delta_p\nabla u(s)\|_{L^2}\|\Delta_{p+j}\nabla u(s)\|_{L^2}\\\times\left\{\sum_{q=3}^{p+2}e^{-2^{2(q-1)}(t-s)}2^{2q}2^{2{(q-p)}}\right\}ds
\ea
$$
Integrating in $t$ and changing order of integration we obtain
$$
\ba
\int_0^T\|C(t)\|_{L^{\infty}}dt\le \\
c\sum_{j=-2}^2\sum_{p=1}^{\infty}\int_0^T\|\Delta_p(\nabla u(s))\|_{L^2}\|\Delta_{p+j}(\nabla u(s))\|_{L^2}\sum_{q=3}^{p+2}2^{2(q-p)}ds \\
\le c\sum_{j=-2}^2\int_0^T\sum_{p=1}^{\infty}\|\Delta_p(\nabla u(s))\|_{L^2}\|\Delta_{p+j}(\nabla u(s))\|_{L^2}ds.
\ea
$$
We have obtained thus
\be
\int_0^T \|C(t)\|_{L^{\infty}}dt \le c E_1.\la{cbound}
\ee

We turn now to the term $I(t)$ of (\ref{it}). This term is made up of two sums,
\be
\ba
I(t) = I_1(t) + I_2(t)\\ 
I_1(t) = \sum_{q\ge 3}\int_0^te^{(t-s)\Delta}\Delta{\mathbb{H}}(D)\Delta_q\left (\sum_{p\ge p'+3}\Delta_p(u(s))\otimes\Delta_{p'}(u(s))\right)ds\\
I_2(t) = \sum_{q\ge 3}\int_0^te^{(t-s)\Delta}\Delta{\mathbb{H}}(D)\Delta_q\left (\sum_{p'\ge p+3}\Delta_p(u(s))\otimes\Delta_{p'}(u(s))\right)ds
\ea
\la{is}
\ee
We will treat $I_1$ because the treatment of $I_2$ is the same, mutatis mutandis. Because $p\ge p'+3$ it follows that, in order to have a nonzero contribution at $q$, the index $p$ must belong to $[q-2,q+2]$, i.e., $p = q +j$ with $j\in [-2,2]$. Then we can write
\be
I_1(t) = \sum_{j=-2}^2\sum_{q\ge 3}\int_0^t e^{(t-s)\Delta}\Delta {\mathbb {H}}(D)\Delta_q (J_q(s))ds\la{ione}
\ee
with
\be
J_q(s) = \Delta_{q+j}(u(s))\otimes S_{q+j-3}(u(s)).
\la{jq}
\ee
For $q\le M$ we estimate
\be
\ba
\|\Delta_q(J_q(s))\|_{L^{\infty}} \le c\|S_{q+j-3}(u(s))\|_{L^{\infty}}\|\Delta_{q+j}(u(s))\|_{L^{\infty}}\\
\le c\left [\|u(s)\|_{L^2} + \sqrt{M+2}\|\nabla u(s)\|_{L^2}\right ]\|\Delta_{q+j}(u(s))\|_{L^{\infty}}
\ea
\la{jl}
\ee
where we used (\ref{sqb}):
\be
\|S_{q+j}(u(s))\|_{L^{\infty}}\le c \left (\|u(s)\|_{L^2} + \sqrt{q+j}\|\nabla u(s)\|_{L^2}\right ).
\la{sqjin}
\ee
Using Bernstein's inequality $\|\Delta_{q+j}u(s)\|_{L^{\infty}}\le c \|\Delta_{q+j}\nabla u(s)\|_{L^2}$ we obtain
\be
\|\Delta_q(J_q(s))\|_{L^{\infty}}\le c\left (\|u(s)\|_{L^2} + \sqrt{M+2}\|\nabla u(s)\|_{L^2}\right)\|\Delta_{q+j}\nabla(u(s))\|_{L^2}
\la{lj}
\ee
For $q\ge M$ we estimate 
\be
\|\Delta_q(J_q(s))\|_{L^{\infty}}\le c\left (\|u(s)\|_{L^2} + \|\Delta u(s)\|_{L^2}\right)2^{-q}\|\Delta_{q+j}\Delta(u(s))\|_{L^2}
\la{hj}
\ee
We write
$$
\int_0^T \|I_1(t)\|_{L^{\infty}}dt \le A +B
$$
with 
$$
A = c\sum_{j=-2}^2\sum_{q=3}^M\int_0^T\|\Delta_q(J_q(s))\|_{L^{\infty}} \left (\int_s^T2^{2q}e^{-(t-s)2^{2(q-1)}}dt \right )ds
$$
and
$$
B = c\sum_{j=-2}^2\sum_{q=M}^{\infty}\int_0^T\|\Delta_q(J_q(s))\|_{L^{\infty}} \left (\int_s^T2^{2q}e^{-(t-s)2^{2(q-1)}}dt \right )ds
$$
We use (\ref{lj}) for $A$:
$$
\ba
A \le c \int_0^T\left (\|u(s)\|_{L^2} + \sqrt{M+2}\|\nabla u(s)\|_{L^2}\right)\left (\sum_{q=3}^M\|\Delta_{q+j}\nabla u(s)\|_{L^2}\right)ds\\\le c\int_0^T \left (\|u(s)\|_{L^2} + \sqrt{M+2}\|\nabla u(s)\|_{L^2}\right)\left( \sqrt{M}\|\nabla u(s)\|_{L^2}\right)ds.
\ea
$$
We have therefore
\be
A\le c M E_1\la{ab}
\ee

We use (\ref{hj}) for $B$:
$$
B \le c\int_0^T\left (\|u(s)\|_{L^2} + \|\Delta u(s)\|_{L^2}\right)\sum_{q\ge M}\left (2^{-q}\|\Delta_{q+j}\Delta(u(s))\|_{L^2}\right)ds
$$
and therefore
$$
B \le c 2^{-M}\int_0^T\|u(s)\|^2_{H^2}ds = c2^{-M}E_2
$$
where
\be
E_2 = \int_0^T\|u(t)\|_{H^2}^2dt.
\la{etwo}
\ee
In view of (\ref{htwo})
\be
B\le c2^{-M}(1+T)\left\{\|\sigma\|_{L^2(0,T; H^1)}^2 + \|u(0)\|^2_{H^1}\right\}
= c2^{-M}(1+T)R_1\la{bb}
\ee
For any $\epsilon>0$, we choose
$$
M =\log_{*}\left (\frac{c(1+T)R_1}{\epsilon}\right )
$$
with $R_1$ given above and we obtain from (\ref{ab}) and  (\ref{bb})
\be
\int_0^T\|I(t)\|_{L^{\infty}}dt \le cE_1\log_{*}\left (\frac{c(1+T)R_1}{\epsilon}\right) + \epsilon\la{ibound}
\ee
The sum of the inequalities (\ref{fbound}), (\ref{s2bound}), 
(\ref{cbound}), (\ref{ibound}) and a straighforward estimate for the linear term carrying the initial data give the inequality (\ref{gineq}) of the theorem.

\section{Coupled Nonlinear Fokker-Planck and \\
Navier-Stokes Systems in 2D}
We consider now the coupling between fluid and particles. The evolution of the density $f$ is governed by a nonlinear
Fokker-Planck equation 
\be
\partial_t f + {\bar{v}}\cdot\nax f + \d(Gf) = \frac{1}{\tau}\D f 
\la{feq} 
\ee 
The coefficient $\tau> 0 $ is the time scale associated with the particles. The microscopic variables $m$ are non-dimensional. 
The tensor $G$ is made of two parts,
\be 
G =  \nag U + W. \la{G} 
\ee 
The $(0,1)$ tensor field $W$ is obtained from the macroscopic gradient of velocity in a linear smooth fashion, given locally as 
\be 
W(x,m,t) = \left
(W_{\alpha}(x,m,t)\right )_{\alpha = 1,\dots, d} = \left (\sum_{i,j
=1}^n c_{\alpha}^{ij}(m)\frac{\partial {\bar {v}}_{i}}{\partial
x_j}(x,t)\right )_{\alpha = 1, \dots, d.} \la{Walphav} 
\ee 
The smooth
coefficients $c^{ij}_{\alpha}(m)$ do not depend on the solution,
time or $x$ and, like the coefficients $\gamma^{(k)}_{ij}$, they are
a constitutive part of the model. 
The potential $U$ is given by 
\be U(x,m,t) =
\frac{b}{\tau}\left({\mathcal K}f\right )(x,m,t) \la{Uv} 
\ee 
where $b$ is a nondimensional measure of the intensity of the inter-particles interaction.  
The nonlocal microscopic interaction potential 
\be 
\left({\mathcal K}f\right)(x,m,t) = \int\limits_{M}K(m,q)f(x,q,t)dq \la{pot} 
\ee is
given by an integral operator with kernel $K(m,q)$ which is a
smooth, time independent, $x$ independent, symmetric function $K:
M\times M \to {\Rm }$ (\cite{ons}). 
The Navier-Stokes equations are
\be
\partial_t v + v\cdot\nax v  + \nax p = \nu\Dx v + \nax\cdot \tau_p\la{nseq}
\ee
with $\nax\cdot v =0$. The added stresses are given by
the relations (\ref{expa}), with (\ref{sigmak}, \ref{sigmasum}). The added stresses are proportional to $kT$ (where $k$ is Botzmann's constant and $T$ is temperature) and have units of energy per unit mass. The density of the fluid is normalized to one.  The particles are advected by
\be
{\bar{v}}(x,t) = \frac{1}{\tau}\int_{(t-\tau)_+}^t v(x,s)ds\la{barv}
\ee
We rescale the Navier-Stokes equations using the length scale $\lambda = \sqrt{\frac{\nu\tau}{\delta}}$ and the time scale $\lambda^2\nu^{-1} = \frac{\tau}{\delta}$, where $\delta$ is the Deborah number, the ratio between the relaxation time scale of the particles and the macroscopic (observation) advective time. We set
$$
v(x,t) = \delta\frac{\lambda}{\tau}u\left (\frac{x}{\lambda}, \frac{t\delta}{\tau}\right)
$$
and we arrive
at (\ref{eq:ns}) with
\be
\sigma = \frac{\tau}{\delta\nu}\tau_p\la{sigmatau}.
\ee
The Fokker-Planck equation becomes
\be
\left(\partial_t + {\bar{u}}\cdot\nax \right)f +\d(Gf) = \frac{1}{\delta}\D f
\la{fe}
\ee
with 
\be
\bar{u}(x,t) = \frac{1}{\delta}\int_{(t-\delta)_+}^t u(x,s)ds,\la{baru}
\ee
and (\ref{G}) with
\be
W(x,m,t) = \left
(W_{\alpha}(x,m,t)\right )_{\alpha = 1,\dots, d} = \left (\sum_{i,j
=1}^n c_{\alpha}^{ij}(m)\frac{\partial {\bar {u}}_{i}}{\partial
x_j}(x,t)\right )_{\alpha = 1, \dots, d.} \la{Walpha}
\ee
and
\be U(x,m,t) =
\frac{b}{\delta}\left({\mathcal K}f\right )(x,m,t) \la{U}
\ee

The forces applied by the particles are
obtained after $f$ is integrated along with smooth coefficients
$\gamma^{(k)}_{ij}$ on $M$ in order to produce $\sigma$. 
Therefore, only very weak regularity of $f$ with respect to the microscopic
variables $m$ is sufficient to control $\sigma$. In order to take
advantage of this, we consider the $L^2(M)$ selfadjoint
pseudodifferential operator \be R = \left (-\D +{\mathbf
I}\right)^{-\frac{s}{2}}\la{R} \ee with $s> \frac{d}{2} + 1$. We
will use the following properties of $R$: \be [R, \nax] =
0,\la{rone} \ee \be R\nag : L^{1}(M)\to L^2(M) \quad {\mbox{is
bounded, }}\la{rtwo} \ee \be R\nag : L^2(M)\to L^{\infty}(M) \quad
{\mbox{is bounded,}} \la{rthree} \ee \be [\nag c, R^{-1}]:
H^{s}(M)\to L^{2}(M) \quad {\mbox{is bounded,}} \la{rfour} \ee for
any smooth function $c: M\to {\mathbb R}$, and \be R: L^{2}(M) \to
H^s(M) \quad {\mbox{is bounded.}} \la{rfive} \ee We differentiate
(\ref{fe}) with respect to $x$, apply $R$, multiply by $R\nax f$ and
integrate on $M$. Let us denote by 
\be N(x,t)^2 = \int\limits_{M}\left
| R\nax f (x,m,t)\right|^2dm \la{Nx} \ee 
the square of the $L^2$
norm of $R\nax f$ on $M$. The following lemma was proved in (\cite{c}):
\begin{lemma}
Let ${\bar{u}}(x,t)$ be a smooth, divergence-free function and let $f$ solve
(\ref{fe}). There exists an absolute  constant $c>0$ (depending only
on dimensions of space, the coefficients $c^{ij}_{\alpha}$ and $M$,
but not on ${\bar{u}}$, $f$, $\delta$) so that

\be 
\left(\partial_t + \bar{u}\cdot \nax \right )N
\le c(|\nax {\bar{u}}| +\frac{1}{\delta})N + c|\nax\nax {\bar{u}}| \la{Nineq} \ee
holds pointwise in $(x,t)$.
\end{lemma}
The proof is given below in the Appendix for completeness. It works
independently of the dimension $n$ of the variables $x$.
The equation obeyed by $\rho =\int_M fdm$ is
\be
(\partial_t + {\bar{u}}\cdot\nax)\rho = 0.
\la{rhoeq}
\ee 
We will take initial densities that 
obey
$$
0\le\rho(x,0)\le 1.
$$
Therefore 
\be
0\le\rho(x,t)\le 1\la{dens}
\ee
continues to be true and, in view of (\ref{sigmalinfty}) and the fact that
$\bar{u}$ is divergence free, it follows from (\ref{rhoeq}) that
\be
\|\sigma(t)\|_{L^r}\le c_r\la{krpb}
\ee
holds if we assume (as we do) that $\rho(x,0)\in L^1({\mathbb R}^2)$.
The inequalities (\ref{sigmalinfty}) and (\ref{naxsigma}) use only the smoothness of the 
coefficients $\gamma$, the relations (\ref{expa}, \ref{sigmak}), the condition
(\ref{sigmasum}) and (\ref{intone}), and therefore they hold throughout the evolution. We know thus that $\sigma$ is bounded in $L^{\infty}$ and that $K_r^{(p)}$ (in ({\ref{krp})) are bounded a priori. (We also know that $\int_0^T\|u(t)\|_{L^{\infty}}dt$ is bounded a priori (see (\ref{almost})). This implies that the support of $\rho$, if initially compact, would expand only a finite amount in finite time. We do not use this for the proof, but obviously this is a physically important a priori quantitative information.)
\begin{thm}
Consider the coupled Fokker-Planck and Navier-Stokes system 
(\ref{feq}), (\ref{nseq}) with arbitrary parameters $\nu, \tau, b >0$. Assume that the initial velocity $v(0)$ is divergence-free and smooth, $v(0)\in W^{2,r}({\mathbb R}^2)\cap L^2({\mathbb R}^2)$ with $r>2$. Assume that the initial distribution of particles $f(x,m,0)$ is non-negative, smooth,  in the sense that
$$
N(x, 0) = \|\nax f(x,\cdot,0)\|_{H^{-s}(M)}\in L^r({\mathbb R}^2)\cap L^2({\mathbb R}^2)
$$ 
for some $s\in {\mathbb {R}}$, and localized, in the sense that
$$
\rho (x,0) = \int_M f(x,m,0)dm
$$
obeys $0\le \rho (x,0)\le 1$ and $\rho(\cdot, 0)\in L^1({\mathbb R}^2)$.
Then the solution of the system (\ref{feq}, \ref{nseq}) exists for all time and is smooth. In particular, the norms of
$$
v\in L^1(dt; W^{1,\infty}({\mathbb R}^2))\cap L^2(dt; W^{2,r}({\mathbb R}^2)),
$$
$$
f\in L^{\infty}(dt; W^{1,\infty}(dx; H^{-s}(M))).
$$
can be bounded a priori in terms of the initial data, for arbitrary large finite intervals of time $[0,T]$.
\end{thm}

\noindent{\bf{Proof.}} Let $T>0$ be given. The sort time existence of solutions and the uniqueness of the solutions can be obtained following classical methods of proof. Let
\be
n(t) = \|\nax \sigma\|_{L^r}^2 + \|\sigma\|^2_{H^1},
\la{nt}
\ee
\be
B(t) = \int_{0}^t\|\nax \sigma\|_{L^r}^2dt + \int_{0}^t\|\sigma\|^2_{H^1}dt,
\la{bt}
\ee
\be
g(t) =\|\nabla {\bar {u}}(t)\|_{L^{\infty}},\la{gt}
\ee
\be
\gamma(t) = \sup_{0\le s\le t}g(s),\la{gamma}
\ee
and
\be
G(t) = \|\nax\nax {\bar{u}}(t)\|_{L^r}^2 + \|{\bar{u}}(t)\|_{H^2}^2.
\la{Gt}
\ee
Using (\ref{Nineq}) we have
\be
\frac{dn}{dt} \le cg(t)n(t) + cG(t)
 \la{nineq}
\ee
Integrating on $(0,t)$
\be
n(t) \le n(0) + c\gamma(t) B(t) + c\int_{0}^t G(s)ds\la{nint}
\ee
In view of (\ref{gineq}) we have
\be
c\gamma(t) \le C_0\left \{1 + \log_{*}(B(t))\right\}
\la{close}
\ee
with $C_0$ a constant depending on $T$ and the initial data.
In view of (\ref{htwo}) and (\ref{grgrr}) we know
\be
\int_{0}^t G(s)ds \le C_1\left\{1 + B(t)\log_{*}B(t) \right\}
\la{Gineq}
\ee
with $C_1$ a constant. Note that
\be
n(t) = \frac{d}{dt}B(t)
\la{dB}
\ee
Thus, from (\ref{nint}), (\ref{close}) and (\ref{Gineq}) we deduce
\be
\frac{d}{dt}B(t) \le  C_2\left\{1  +  B(t)\log_{*}B(t)\right\} 
\la{bla} 
\ee
holds with $C_2$ a constant that depends on $T$ and the values
$\|\nabla\omega(0)\|_{L^r}$, $\|\omega(0)\|_{L^2}$, $\|u(0)\|_{L^2}$ and $n(0)$.
This produces a pointwise-in-time a priori finite bound for $B(t)$ on the interval $[0,T]$,
and retracing our steps, via (\ref{Gineq}) and (\ref{nint}), on $n(t)$. Once the forces in the two-dimensional Navier-Stokes equations are known to be thus bounded, it follows (from (\ref{gineq})--but also much easier, from energy estimates) that the solution of the Navier-Stokes equation is smooths as stated.

\section{Appendix: Proof of Lemma 3}
The evolution equation of $N$ is \be \frac{1}{2}\left (\partial_t +
{\bar{u}}\cdot\nax \right ) N^2  =  -D  + I + II  + III + IV
\la{neq} \ee with \be D= \epsilon\int\limits_{M}\left |
\nag R \nax f \right |^2dm \la{D} \ee 
\be I =  -\frac{\partial
{\bar{u}}_j}{\partial x_k}\int\limits_{M} \left (R\frac{\partial f}{\partial
x_j}\right ) \left (R\frac{\partial f}{\partial x_k}\right )dm
\la{I} \ee 
\be II = -\sum\limits_{\alpha =1}^2 (\nax \frac{\partial
{\bar{u}}_i}{\partial {x_j}}) \int\limits_{M}(R\d(c^{ij}_{\alpha} f))(\nax
Rf) dm, \la{II} \ee \be III = - \sum\limits_{\alpha
=1}^2\frac{\partial {\bar{u}}_i}{\partial x_j} \int\limits_{M}(R \d
(c^{ij}_{\alpha}\nax f))(R\nax f) dm, \la{III} \ee and \be IV =
-\frac{b}{\delta}\int\limits_{M}R\d(\nax \left\{f\nag({\mathcal
K}f)\right\})R\nax f dm. \la{IV} \ee

Now we start estimating these terms. We will use repeatedly (\ref{intone}) and (\ref{dens}). In view of the fact that $D\ge 0$, we may discard this term.  Clearly \be |I|\le c |\nax {\bar{u}} | N^2. \la{Ineq} \ee In
order to bound $II$ we use (\ref{rtwo}) to bound
$$
\| R\nag (c^{ij}_{\alpha} f)\|_{L^2(M)}\le c \|f\|_{L^1(M)} = c
$$
so that we have \be |II|\le c\left| \nax\nax {\bar{u}} \right| N. \la{IIneq}
\ee In order to bound $III$ we need to use the commutator carefully.
We start by writing
$$
R\d(c^{ij}_{\alpha}\nax f) = R\d(c^{ij}_{\alpha}R^{-1}R\nax f) =
$$
$$
\d(c^{ij}_{\alpha} R\nax f) + \left [R\d c^{ij}_{\alpha},
R^{-1}\right ]R\nax f.
$$
The second term obeys
$$
\|\left [R\d c^{ij}_{\alpha}, R^{-1}\right ]R\nax f\|_{L^2(M)} \le c
N
$$
because, in view of (\ref{rfour}) and (\ref{rfive}) one has that
$$
\left [R\d c^{ij}_{\alpha}, R^{-1}\right ]: L^2(M)\to L^2(M)
\quad{\mbox{is bounded.}}
$$
The first term needs to be integrated against $R\nax f$ and
integration by parts gives
$$
\int\limits_{M}(\d (c^{ij}_{\alpha} R\nax f)) R\nax f dm =
\frac{1}{2}\int\limits_{M}(\d c^{ij}_{\alpha})|R\nax f|^2dm.
$$
We obtain thus \be |III| \le c \left |\nax {\bar{u}}\right | N^2 \la{IIIneq}
\ee The term $IV$ is split in two terms, $IV= A+B$ \be A=
-\frac{b}{\delta}\int\limits_{M}R\d(\left\{(\nax f)\nag({\mathcal
K}f)\right\})R\nax f dm \la{A} \ee and \be B=
-\frac{b}{\delta}\int\limits_{M}R\d(\left\{f\nag({\mathcal K}\nax
f)\right\})R\nax f dm. \la{B} \ee 
The $(0,1)$ tensor $\Phi(x,m,t) =
(\nag{\mathcal K}f)(x,m,t)$ is smooth in $m$ for fixed $x,t$ and
$$
\|\Phi (x,\cdot, t)\|_{W^{s,\infty}(M)} \le c_s
$$
holds for any $s$, with $c_s$ depending only on the kernel $K$. We
write the term $A$
$$
\begin{array}{c}
A = -\frac{b}{\delta}\int\limits_{M}R\d(\left\{(\nax f)
\Phi\right\})R\nax fdm  \\
= \frac{b}{\delta}\int\limits_{M}R^{-1}(R\nax f)\{\Phi\cdot\nag
R^2\nax
f))dm \\
= -\frac{b}{2\delta}\int\limits_{M}\d\left\{\Phi\right\}\left |R\nax
f\right |^2dm + \frac{b}{\delta}\int\limits_{M}(R\nax f)\left [
R^{-1}, \Phi\nag \right ]R(R\nax f)dm.
\end{array}
$$
In view of (\ref{rfour}), (\ref{rfive}), the operator
$$
\left [R^{-1}, \Phi\nag \right]R : L^2(M) \to L^2(M)
$$
is bounded with norm bounded by an a priori constant. It follows
that
$$
|A| \le \frac{cb}{\delta}N^2(x,t)
$$
holds. The term $B$ is easier to bound, because
$$
({\mathcal K}\nax f)(x,m,t) = \int\limits_{M}R^{-1}K(m,n) R\nax
f(x,n,t)dn
$$
and thus
$$
\|(\nag{\mathcal K}\nax f)(x,\cdot, t)\|_{L^{\infty}(M)}\le cN(x,t).
$$
Using (\ref{rtwo}) it follows that
$$
|B|\le  \frac{cb}{\delta}N^2(x,t)
$$
and consequently \be |IV| \le \frac{cb}{\delta}N^2(x,t).\la{IVneq} \ee

Putting together the inequalities (\ref{Ineq}), (\ref{IIneq}),
(\ref{IIIneq}) and (\ref{IVneq}) we finished the proof of the lemma.

\vspace{1cm}
\noindent{\bf{Acknowledgment.}}\,\,\, The work of P.C. is partially supported by NSF-DMS grant 0504213.\,\, The work of C.F. is partially supported by NSF-DMS grant 0245242.\,\,\,The work of  E.S.T. is supported in part by the NSF grant no.  
DMS-0504619, the BSF grant no. 2004271, and by the MAOF Fellowship 
of the Israeli Council of Higher Education.


\begin{thebibliography}{99}
\bibitem{brega} H.Br\'{e}zis and T.Gallouet, Nonlinear Schr\"{o}dinger
evolution equations,  Nonlinear Anal.  {\bf 4}  (1980), no. 4, 677--681.
\bibitem{CH98} J. Y. Chemin, {\em Perfect Incompressible Fluids},
Oxford Lecture Series in Mathematics and its Applications, {\bf 14},
 Clarendon Press, Oxford University Press, New York, 1998.
\bibitem{CHMA01} J.Y. Chemin, and N. Masmoudi  About lifespan of
 regular solutions of equations related to viscoelastic fluids,
 SIAM J. Math. Anal.  {\bf 33}  (2001),  no. 1, 84--112.
\bibitem{c} P. Constantin, Nonlinear Fokker-Planck Navier-Stokes
systems, Commun. Math. Sciences, {\bf 3} (4) (2005), 531-544.
\bibitem{cfbook} P. Constantin, C. Foias, {\em Navier-Stokes Equations}, U.
Chicago Press, Chicago 1988.
\bibitem{doied} M. Doi, S.F. Edwards, {\em The Theory of Polymer Dynamics},
Oxford University Press, Oxford 1988.
\bibitem{ez} W. E, T.J. Li, P-W. Zhang, Well-posedness for the dumbell
model of polymeric fluids, Commun. Math. Phys. {\bf 248} (2004),
409-427.
\bibitem{hel} S. Helgason, {\em Differential Geometry, Lie Groups and
Symmetric spaces}, Academic Press, London 1978.
\bibitem{hor} L. H\"{o}rmander, {\em The Analysis of Linear Partial
Differential Operators}, vol. {\bf 3}, Springer-Verlag, Berlin, Heidelberg,
New York, Tokyo, 1985.
\bibitem{jourdain} B. Jourdain, T. Lelievre, C. Le Bris, Numerical
analysis of micro-macro simulations of polymeric flows: a simple
case, Math. Models in Appl. Science {\bf 12} (2002), 1205-1243.
\bibitem{jourdainl} B. Jourdain, T. Lelievre, C. Le Bris, Esistence of
solutions for a micro-macro model of polymeric fluid: the FENE
model, preprint (2005).
\bibitem{lemarie} P.G. Lemari\'{e}-Rieusset, {\em Recent Developments in the
Navier-Stokes Problem}, Chapmann and Hall/CRC, Research Notes in
Mathematics {\bf 431}, CRC, Boca Raton, 2002.

\bibitem{linliuzang} F-H. Lin, C. Liu, P. Zhang, On hydrodynamics of viscoelastic fluids, CPAM {\bf 58}, (2005), 1437-1471.
 
\bibitem{lizz} T.Li, H. Zhang, P-W. Zhang, Local existence for the
dumbell model of polymeric fluids, Commun. PDE {\bf 29} (2004),
903-923.
\bibitem{ons} L. Onsager, The effects of shape on the interaction of
colloidal particles, Ann. N.Y. Acad. Sci {\bf 51} (1949), 627-659.
\bibitem{otto-tzavaras} F. Otto, A. Tzavaras, Continuity of velocity
gradients in suspensions of rod-like molecules, SFB Preprint {\bf
147} (2004).
\bibitem{renardy} M. Renardy, An existence theorem for model equations
resulting from kinetic theories of polymer solutions, SIAM J. Math.
Analysis {\bf 23} (1991), 313-327.
\bibitem{sideris} T. Sideris, B. Thomasses, Global existence for 3D
incompressible isotropic elastodynamics via the incompressible limit, CPAM {\bf 58} (2005), 750-788.

\end{thebibliography}
\end{document}